\documentclass[12pt]{amsart}
 
\pdfoutput=1
\usepackage[margin=1in]{geometry}
\usepackage{amsmath,amsthm,amssymb,amsrefs,bbm,color,esvect,float,graphicx,mathrsfs}
\usepackage{epstopdf}
\AppendGraphicsExtensions{.tif}

\newenvironment{theorem}[2][Theorem]{\begin{trivlist}
\item[\hskip \labelsep {\bfseries #1}\hskip \labelsep {\bfseries #2.}]}{\end{trivlist}}
\newenvironment{definition}[2][Definition]{\begin{trivlist}
\item[\hskip \labelsep {\bfseries #1}\hskip \labelsep {\bfseries #2.}]}{\end{trivlist}}

\newenvironment{corollary}[2][Corollary]{\begin{trivlist}
\item[\hskip \labelsep {\bfseries #1}\hskip \labelsep {\bfseries #2.}]}{\end{trivlist}}
\newenvironment{remark}[2][Remark]{\begin{trivlist}
\item[\hskip \labelsep {\bfseries #1}\hskip \labelsep {\bfseries #2.}]}{\end{trivlist}}
\theoremstyle{definition}

\begin{document}

\title[A limited-range Calder\'on-Zygmund theorem]{A limited-range Calder\'on-Zygmund theorem}
\author{Loukas Grafakos}
\address{Loukas Grafakos, Department of Mathematics, University of Missouri, Columbia, MO 65211, USA}
\email{grafakosl@missouri.edu}
\author{Cody B. Stockdale}
\address{Cody B. Stockdale, Department of Mathematics and Statistics, Washington University in St. Louis, One Brookings Drive, St. Louis, MO, 63130, USA}
\email{codystockdale@wustl.edu}
\thanks{The second author would like to acknowledge the Simons Foundation}

\maketitle

\begin{abstract}
For a limited  range of indices $p$, 
we obtain  $L^p(\mathbb{R}^n)$ boundedness for singular integral operators whose kernels satisfy a condition weaker than the typical H\"ormander smoothness estimate. These operators are assumed to be 
bounded (or weakly bounded) on $L^{s}(\mathbb{R}^n)$ for some index $s$. Our estimates are 
obtained via interpolation from the appropriate weak-type estimates. We provide two proofs of this result. 
One proof  is based on  the Calder\'on-Zygmund decomposition, while the other uses ideas of Nazarov, Treil, and Volberg.
\end{abstract}


\section{Introduction}
The classical theory of singular integral operators was introduced by Calder\'on and Zygmund in \cite{CZ1952} and says that for certain kernels defined on $\mathbb{R}^n\setminus \{0\}$, the weak-type $(1,1)$ bound holds for the associated singular integral operator, assuming that an $L^s(\mathbb{R}^n)$ bound is known for some $1<s\leq\infty$. H\"ormander extended this theory in \cite{H1960} to more general kernels $K$ satisfying the smoothness condition $$[K]_{H}:= \sup_{y\in \mathbb{R}^n} \int_{|x|\ge 2|y|} |K(x-y)-K(x)|\, dx <\infty.$$ The H\"ormander condition is an $L^1(\mathbb{R}^n)$-type smoothness condition and has some variants. For example, Watson introduced the following $L^r(\mathbb{R}^n)$ versions in \cite{W1990}: for $1\le r \le \infty$, we say a kernel $K$ is in the class $H^r$ if
$$[K]_{H^r}:=\sup_{R>0}\sup_{\substack{ y\in \mathbb{R}^n \\ |y|\leq R}}  \sum_{m=1}^{\infty}  (2^mR)^{\frac{n}{r'}}
\Bigg[\int\limits_{\substack{|x|\ge 2^m R    \\ |x|<2^{m+1} R}} |K(x-y)-K(x)|^r  dx \Bigg]^{\frac{1}{r}} <\infty,$$ where $r'$ is the H\"older conjugate of $r$. Observe that Watson's condition coincides with H\"ormander's condition when $r=1$, and for $r_1,r_2 \in [1,\infty]$ with $r_1\leq r_2$, $$H^{r_2} \subseteq H^{r_1}  \subseteq  H^1=H.$$ 

In this paper, we focus on a different set of $L^r(\mathbb{R}^n)$-adapted conditions defined as follows.
\begin{definition}{1} 
Let $1\leq r \leq \infty$. A kernel $K$ defined on $\mathbb{R}^n\setminus \{0\}$ is in the class $H_r$ if
$$[K]_{H_r}:=\sup_{R>0} \bigg[\frac{1}{v_nR^n}\int_{|y|\le R} \bigg(
 \int_{|x|\ge 2 R}   |K(x-y)-K(x)|   \, dx\bigg)^{r} dy \bigg]^{\frac{1}{r}} <\infty,$$ where $v_n$ is the volume of the unit ball $B(0,1)$ in $\mathbb{R}^n$.
\end{definition}
Notice that this condition coincides with the H\"ormander condition when $r=\infty$. Moreover, for $r_1,r_2 \in [1,\infty]$ with $r_1\leq r_2$, $$H =H_{\infty} \subseteq  H_{r_2}\subseteq H_{r_1},$$ meaning the $H_r$ conditions are weaker than H\"ormander's smoothness condition.

We prove boundedness results for the associated singular integral operators.
\begin{definition}{2}
Let $K \in H_{r}$ for some $1\leq r \leq \infty$ and suppose $K$ satisfies the size estimate $|K(x)|\leq \frac{A}{|x|^n}$ for all $x\neq 0$. We associate $K$ with a linear operator $T$ given by $$Tf(x) = \int_{\mathbb{R}^n} K(x-y) f(y)\, dy$$ for smooth functions $f$ and $x \not \in \text{supp} f$.
\end{definition}
Notice that this definition also makes sense if $f$ is an integrable, compactly supported function and $x \not \in \text{supp} f$. Moreover, there is no unique way to define $Tf$ in terms of $K$ for general functions $f$ (see the relevant discussions in \cites{Grafakos1, Grafakos2,Stein}).

If $K \in H=H_{\infty}$, H\"ormander proved that given $1<s\leq \infty$, $L^s(\mathbb{R}^n)$ bounds for $T$ imply the weak-type $(1,1)$ bound, and hence $L^p(\mathbb{R}^n)$ bounds for all $1< p< \infty$. In this note, we prove the following variant of this result, where weak-type $(1,1)$ is replaced by weak-type $(q,q)$. 
\begin{theorem}{1}
Let $1 \leq q < \infty$, $K \in H_{q'}$, and $|K(x)|\leq \frac{A}{|x|^{n}}$ for all $x\neq 0$. If the associated singular integral operator $T$ is bounded on $L^s(\mathbb{R}^n)$ for some $s \in (q,\infty]$ with bound $B$, then $T$ maps $L^q(\mathbb{R}^n)$ to $L^{q,\infty}(\mathbb{R}^n)$ with bound at most a constant multiple of $B+[K]_{H_{q'}}$. That is, $$\|Tf\|_{L^{q,\infty}(\mathbb{R}^n)}:=\sup_{\alpha>0}\alpha|\{|Tf|>\alpha\}|^{\frac{1}{q}}\leq C_{n,s,q}(B+[K]_{H_{q'}})\|f\|_{L^q(\mathbb{R}^n)}$$ for all $f \in L^q(\mathbb{R}^n)$.
\end{theorem}

We give two proofs of Theorem 1. The first proof uses the $L^q(\mathbb{R}^n)$ version of the Calder\'on-Zygmund decomposition and is an adaptation of the classical proof given in \cite{CZ1952}. The second proof is motivated by Nazarov, Treil, and Volberg's proof for the weak-type $(1,1)$ inequality in the nonhomogeneous setting, given in \cite{NTV1998}. Adaptations of the proof in the nonhomogeneous setting are needed in our setting; some modifications include ideas that can be found in \cite{Stein}. See \cites{S2018,S2019, SW2019} for other applications of the Nazarov, Treil, and Volberg technique to multilinear and weighted settings. Refer to \cites{GT2002,LOPTTG2009,OPR2016} for related results regarding multilinear and weighted Calder\'on-Zygmund theory.

By interpolation we obtain the following corollary. 
\begin{corollary}{1}
Under the hypotheses of Theorem 1, the operator $T$ is bounded
on $L^p(\mathbb{R}^n)$ for $p$ in the interval $\big(\min (s',q), \max (q',s)  \big)$. 
\end{corollary}

\begin{remark}{1}
The constant $A$ does not appear in the conclusion of Theorem 1. The estimate $|K(x)|\le \frac{A}{|x|^n}$ is only needed to ensure that the operator $T$ is well-defined for a dense class of functions.
\end{remark}

 If $q>1$  and $s<\infty$, then the interval $\big(\min (s',q), \max (q',s)  \big)$ is properly contained in $(1,\infty)$.  Hence in this case, we obtain $L^p(\mathbb{R}^n)$ estimates for a limited-range of 
 values of $p$.  Prior to this work, 
 other ``limited-range'' versions of the Calder\'on-Zygmund theorem  appeared in
   Baernstein and Sawyer \cite{BS1985},
Carbery \cite{C1986},
Seeger \cite{S1988}, and 
Grafakos, Honz\'ik, Ryabogin \cite{GHR}. . 


\section{Calder\'on-Zygmund Decomposition Method}
The first proof of Theorem 1 relies on the $L^q(\mathbb{R}^n)$ version of the Calder\'on-Zygmund decomposition. See \cites{Grafakos1,Grafakos2,Stein} for details on the decomposition.

\begin{proof}
Fix $f \in L^q(\mathbb{R}^n)$ and $\alpha>0$. We will show that $$|\{|Tf|>\alpha\}| \leq C_{n,s,q}(B+[K]_{H^{q'}})^q\alpha^{-q}\|f\|_{L^q(\mathbb{R}^n)}^q.$$ Apply the $L^q(\mathbb{R}^n)$-form of the Calder\'on-Zygmund decomposition to $f$ at height $\gamma\alpha$ (the constant $\gamma>0$ will be chosen later), to write $f=g+b=g+\sum_{j=1}^{\infty}b_j$, where
\begin{enumerate}
\addtolength{\itemsep}{0.2cm}
\item[(1)] $\, \|g \|_{L^{\infty}(\mathbb{R}^n)} \le 2^{\frac{n}{q}} \gamma\alpha$ and $\|g \|_{L^q(\mathbb{R}^n)}\le   \|f \|_{L^q(\mathbb{R}^n)}$,
\item[(2)] \, the $b_j$ are supported on pairwise disjoint cubes $Q_j$ satisfying $\sum_{j=1}^{\infty} |Q_j| \leq (\gamma\alpha)^{-q}  \|f \|_{L^q(\mathbb{R}^n)}^q$,
\item[(3)] $\,  \|b_j \|_{L^q(\mathbb{R}^n)}^q \le 2^{n+q} (\gamma\alpha)^q |Q_j|$,
\item[(4)] $\, \int_{Q_j} b_j(x)\, dx= 0$, and
\item[(5)]  $\,  \|b \|_{L^q(\mathbb{R}^n)} \leq 2^{\frac{n+q}{q}}
  \|f \|_{L^q(\mathbb{R}^n)} $ and  $ \|b \|_{L^1(\mathbb{R}^n)} \le 2 
(\gamma\alpha)^{1-q} \|f \|_{L^q(\mathbb{R}^n)}^q$.
\end{enumerate}
Now, $$|\{|Tf|>\alpha\}|\leq \left|\left\{|Tg|>\frac{\alpha}{2}\right\}\right|+\left|\left\{|Tb|>\frac{\alpha}{2}\right\}\right|.$$

Assume first that $s<\infty$. Choose $\gamma=(B+[K]_{H_{q'}})^{-1}$. Using Chebyshev's inequality, the bound of $T$ on $L^s(\mathbb{R}^n)$, property (1), and trivial estimates, we have that 
\begin{align*}
\left| \left\{ |Tg|>\frac{\alpha}{2}\right\}\right| & \leq 2^s\alpha^{-s}\|Tg\|_{L^s(\mathbb{R}^n)}^s  \\
&\leq (2B)^s\alpha^{-s} \| g\|_{L^s(\mathbb{R}^n)}^s  \\
& \leq 2^{s-n+\frac{ns}{q}}B^s\alpha^{-s}  
 (\gamma\alpha)^{s-q}\| g \|_{L^q(\mathbb{R}^n)} ^q  \\
 & \leq  2^{s-n +\frac{ns}{q}}(B+[K]_{H_{r'}})^q \alpha^{-q} \|f\|_{L^q(\mathbb{R}^n)}^q  .
\end{align*}

We next control the second term. Let $c_j$ denote the center of $Q_j$, let $Q_j^*:=Q(c_j,2\sqrt{n}l(Q_j))$ be the cube centered at $c_j$ and having side length $2\sqrt{n}$ times the side length of $Q_j$, and set $\Omega^*:=\bigcup_{j=1}^{\infty}Q_j^*$. Then $$\left|\left\{|Tb|>\frac{\alpha}{2}\right\}\right|\leq |\Omega^*|+\left|\left\{x \in \mathbb{R}^n\setminus \Omega^*: |Tb(x)|>\frac{\alpha}{2}\right\}\right|.$$ Notice that since $|Q_j^*|=(2\sqrt{n})^n|Q_j|$ and by property (2), we have $$|\Omega^*|\leq\sum_{j=1}^{\infty}|Q_j^*|= (2\sqrt{n})^n\sum_{j=1}^{\infty}|Q_j|\leq(2\sqrt{n})^n(B+[K]_{H_q'})^q\alpha^{-q}\|f\|_{L^q(\mathbb{R}^n)}^q.$$ 

It remains to control the last term. Use Chebyshev's inequality, property (4), Fubini's theorem, H\"older's inequality, property (3), and property (2) to estimate
\begin{align*}
&  \left|\left\{\mathbb{R}^n\setminus \Omega^*: |Tb|>\frac{\alpha}{2}\right\}\right|   \leq 2\alpha^{-1} \int_{\mathbb{R}^n\setminus\Omega^*}
| Tb(x)|\, dx \\
& \quad\quad\leq 2\alpha^{-1}\sum_{j=1}^{\infty} \int_{\mathbb{R}^n\setminus\Omega^*}|Tb_j(x)| dx \\
& \quad\quad\leq 2\alpha^{-1} \sum_{j=1}^{\infty}  \int_{Q_j} 
 \bigg[\int_{\mathbb{R}^n\setminus\Omega^*} |K(x-y)-K(x-c_{j})| \,  dx\bigg] |b_j(y)|\, dy \\
& \quad\quad\le 2\alpha^{-1} \sum_{j=1}^{\infty} 
 \bigg\| \int_{\mathbb{R}^n\setminus\Omega^*} |K(x-\cdot)-K(x-c_{j})| \,  dx\bigg\|_{L^{q'}(Q_j)} \| b_j \|_{L^q} \\ 
 & \quad\quad\le 2\alpha^{-1} \sup_{j\in\mathbb{N}} \bigg\| \int_{\mathbb{R}^n\setminus\Omega^*} |K(x-\cdot)-K(x-c_{j})| \,  dx\bigg\|_{L^{q'}\left(Q_j, \frac{dy}{|Q_j|} \right)} \sum_{j=1}^{\infty}  |Q_j|^{\frac{1}{q'}}
 \| b_j \|_{L^q} \\ 
&\quad\quad \leq 2^{\frac{n}{q}+2}\gamma 
\sup_{j\in\mathbb{N}} \bigg\| \int_{\mathbb{R}^n\setminus\Omega^*} |K(x-\cdot)-K(x-c_{j})| \,  dx\bigg\|_{L^{q'}\left(Q_j, \frac{dy}{|Q_j|} \right)}\sum_{j=1}^{\infty}|Q_j|\\
&\quad\quad \leq2^{\frac{n}{q}+2}\gamma^{1-q}\alpha^{-q}\|f\|_{L^q(\mathbb{R}^n)}^q\sup_{j\in\mathbb{N}} \bigg\| \int_{\mathbb{R}^n\setminus\Omega^*} |K(x-\cdot)-K(x-c_{j})| \,  dx\bigg\|_{L^{q'}\left(Q_j, \frac{dy}{|Q_j|} \right)}.
\end{align*}
For each $j$, setting $R_j=\frac{\sqrt{n}}{2}l(Q_j)$, we have
$$
Q_j \subseteq  B(c_{j},  R_j) \subseteq  B(c_{j}, 2R_j)  \subseteq  Q_j^*,
$$ where $B(x,r)$ denotes the ball centered at $x$ and with radius $r$. Then the factor involving the supremum is less than or equal to
 $$\sup_{j\in\mathbb{N}}\bigg[ \int_{B(c_{j}, R_j)} 
  \bigg(\int_{\mathbb{R}^n\setminus B(c_{j}, 2R_j)}  |K(x-y)-K(x-c_{j})|   dx
 \bigg)^{q'}\!\! \frac{dy}{|Q_j|} \bigg]^{\frac{1}{q'}},$$
which is bounded by $\left(\frac{\sqrt{n}}{2}\right)^nv_n [K]_{H_{q'}}$ by changing variables $x'=x-c_{j}$, 
$y'=y-c_{j}$ and 
by replacing the supremum over $R_j $ by the supremum over all $R>0$. 

Putting all of the estimates together, we get $$|\{|Tf|>\alpha\}|\leq \left(2^{s-n +\frac{ns}{q}}+(2\sqrt{n})^n+2^{\frac{n}{q}+2-n}n^{\frac{n}{2}}\right)(B+[K]_{H_{q'}})^q\alpha^{-q}\|f\|_{L^q(\mathbb{R}^n)}^q.$$

When $s=\infty$, set $\gamma=  2^{-\frac{n}{q}} (4([K]_{H_{q'}}+B))^{-1}$. Then 
$$
\|Tg\|_{L^{\infty}(\mathbb{R}^n)} \leq B \|g\|_{L^{\infty}(\mathbb{R}^n)} \leq 2^{\frac{n}{q}}B\gamma\alpha\leq \frac{\alpha}{4},
$$
so $$\left|\left\{|Tg|>\frac{\alpha}{2}\right\}\right|=0.$$ The part of the argument involving $\left\{|Tb|>\frac{\alpha}{2}\right\}$ is the same as in the case $s<\infty$. 
\end{proof}


\section{Method of Nazarov, Treil, and Volberg }
We provide a second proof of Theorem 1. This proof is motivated by the argument given by Nazarov, Treil, and Volberg in \cite{NTV1998}. See also \cites{S2018,S2019,SW2019} for other applications of this technique.

\begin{proof}
Fix $f\in L^q(\mathbb{R}^n)$ and $\alpha>0$. We will show that $$|\{|Tf|>\alpha\}|\leq C_{n,s,q}(B+[K]_{H_{q'}})^q\alpha^{-q}\|f\|_{L^q(\mathbb{R}^n)}^q.$$ By density, we may assume $f$ is a nonnegative continuous function with compact support. Set $$\Omega:=\left\{M(f^q)>(\gamma \alpha)^q\right\}$$ where $\gamma>0$ is to be chosen later and where $M$ denotes the Hardy-Littlewood maximal operator. Apply a Whitney decomposition to write $$\Omega=\bigcup_{j=1}^{\infty}Q_j,$$ a disjoint union of dyadic cubes 
where $$2\text{diam}(Q_j)\leq d(Q_j,\mathbb{R}^n\setminus \Omega)\leq 8\text{diam}(Q_j).$$ Put $$g:=f\mathbbm{1}_{\mathbb{R}^n\setminus \Omega}, \quad\quad b:=f\mathbbm{1}_{\Omega}, \quad\quad \text{and} \quad\quad b_j:=f\mathbbm{1}_{Q_j}.$$ Then $$f=g+b=g+\sum_{j=1}^{\infty}b_j,$$ where we claim that 
\begin{enumerate}
\addtolength{\itemsep}{0.2cm}
\item[(1)] $\, \|g\|_{L^{\infty}(\mathbb{R}^n)}\leq \gamma \alpha$ and  $\|g\|_{L^q(\mathbb{R}^n)}\leq \|f\|_{L^q(\mathbb{R}^n)}$,
\item[(2)] \, the $b_j$ are supported on pairwise disjoint cubes $Q_j$ satisfying 
$$
\sum_{j=1}^{\infty}|Q_j|\leq 3^n(\gamma \alpha)^{-q}\|f\|_{L^q(\mathbb{R}^n)}^q ,
$$
\item[(3)] $\|b_j\|_{L^q(\mathbb{R}^n)}^q \leq (17\sqrt{n})^n(\gamma \alpha)^q |Q_j|$, and
\item[(4)] $\|b\|_{L^q(\mathbb{R}^n)}\leq\|f\|_{L^q(\mathbb{R}^n)}$ and $\|b\|_{L^1(\mathbb{R}^n)}\leq (17\sqrt{n})^{\frac{n}{q}}3^n(\gamma\alpha)^{1-q}\|f\|_{L^q(\mathbb{R}^n)}^q$.
\end{enumerate}
Indeed, since for any $x \not \in \Omega$, we have $$|g(x)|^q=|f(x)|^q\leq M(f^q)(x) \leq (\gamma\alpha)^q,$$  it follows that $\|g\|_{L^{\infty}(\mathbb{R}^n)}\leq \gamma\alpha$. Since $g$ is a restriction of $f$, we have $\|g\|_{L^q(\mathbb{R}^n)}\leq\|f\|_{L^q(\mathbb{R}^n)}$, and so (1) holds. Using the weak-type $(1,1)$ bound for $M$ with $\|M\|_{L^1(\mathbb{R}^n)\rightarrow L^{1,\infty}(\mathbb{R}^n)}\leq 3^n$, we obtain 
property (2) as follows $$\sum_{j=1}^{\infty}|Q_j|=|\Omega|\leq 3^n(\gamma\alpha)^{-q}\|f\|_{L^q(\mathbb{R}^n)}^q.$$  

Addressing (3) and (4), let $Q_j^*:=Q(c_j,17\sqrt{n}l(Q_j))$ be the cube with the same center as $Q_j$ but side length $17\sqrt{n}$ times as large. Then $Q_j^*\cap (\mathbb{R}^n\setminus\Omega)\neq \emptyset$, so there is a point $x \in Q_j^*$ such that $M(f^q)(x)\leq (\gamma\alpha)^q$. In particular, $\int_{Q_j^*}|f(y)|^qdy\leq (\gamma\alpha)^q|Q_j^*|$. Since $|Q_j^*|=(17\sqrt{n})^n|Q_j|$, we have $$\|b_j\|_{L^q(\mathbb{R}^n)}^q=\int_{Q_j}|f(y)|^qdy\leq \int_{Q_j^*}|f(y)|^qdy\leq(\gamma\alpha)^q|Q_j^*|=(17\sqrt{n})^n(\gamma\alpha)^q|Q_j|.$$ This proves (3). We use H\"older's inequality, property (3), and property and (2) to justify property (4) $$\|b\|_{L^1(\mathbb{R}^n)}=\sum_{j=1}^{\infty}\|b_j\|_{L^1(\mathbb{R}^n)}\leq\sum_{j=1}^{\infty}\|b_j\|_{L^q(\mathbb{R}^n)}|Q_j|^{\frac{1}{q'}}\leq (17\sqrt{n})^{\frac{n}{q}}(\gamma\alpha)\sum_{j=1}^{\infty}|Q_j|$$ $$\leq (17\sqrt{n})^{\frac{n}{q}}3^n(\gamma\alpha)^{1-q}\|f\|_{L^q(\mathbb{R}^n)}^q.$$

Now,
$$|\{|Tf|>\alpha\}|\leq \left|\left\{|Tg|>\frac{\alpha}{2}\right\}\right|+\left|\left\{|Tb|>\frac{\alpha}{2}\right\}\right|.$$ Assume first that $s<\infty$. Choose $\gamma=(B+[K]_{H_{q'}})^{-1}$. Use Chebyshev's inequality, the bound of $T$ on $L^s(\mathbb{R}^n)$, and property (1) to see 
\begin{align*}
\left|\left\{|Tg|>\frac{\alpha}{2}\right\}\right|&\leq 2^s\alpha^{-s}\|Tg\|_{L^s(\mathbb{R}^n)}^s\\
&\leq (2B)^s\alpha^{-s}\|g\|_{L^s(\mathbb{R}^n)}^s\\
&\leq (2B)^s(\gamma\alpha)^{s-q}\alpha^{-s}\|g\|_{L^q(\mathbb{R}^n)}^q\\
&\leq 2^s(B+[K]_{H_{q'}})^q\alpha^{-q}\|f\|_{L^q(\mathbb{R}^n)}^q.
\end{align*}

We will now control the second term. Let $E_j$ be a concentric dilate of $Q_j$; precisely,
$$E_{j}:=Q(c_{j},r_{j}),$$ where $c_j$ is the center of $Q_j$ and $r_{j}>0$ is chosen so that $|E_{j}|=\frac{1}{(17\sqrt{n})^{\frac{n}{q}}\gamma\alpha}\int_{Q_j}b_j(x)\,dx$. Note that such $E_j$ exist since the function $r \mapsto |Q(x,r)|$ is continuous for each $x \in \mathbb{R}^n$. Applying H\"older's inequality and property (3), we have $$|E_j|=\frac{1}{(17\sqrt{n})^{\frac{n}{q}}\gamma\alpha}\int_{Q_j}b_j(x)\,dx\leq \frac{1}{(17\sqrt{n})^{\frac{n}{q}}\gamma\alpha}|Q_j|^{\frac{1}{q'}}\|b_j\|_{L^q(\mathbb{R}^n)} \leq |Q_j|.$$ Since $E_j$ is a cube with the same center as $Q_j$ and since $|E_j|\leq |Q_j|$, the containment $E_j\subseteq Q_j$ holds. In particular, the $E_j$ are pairwise disjoint. Set $$E:=\bigcup_{j=1}^{\infty}E_{j}.$$ Then
\[
\left|\left\{|Tb|>\frac{\alpha}{2}\right\}\right| 
\le \text{I}+\text{II}+\text{III},
\]
where 
\begin{align*}
\text{I}& = |\Omega |, \\
\text{II}& = \left|\left\{x\in\mathbb{R}^n\setminus \Omega:\left|T\left(b-(17\sqrt{n})^{\frac{n}{q}}\gamma\alpha\mathbbm{1}_E\right)(x)\right|>\frac{\alpha}{4}\right\}\right|, \,\, \text{and}\\
\text{III}& =\left|\left\{ (17\sqrt{n})^{\frac{n}{q}}\gamma\alpha|T(\mathbbm{1}_E)|>\frac{\alpha}{4}\right\}\right|.
\end{align*}

The control of I follows from property (2), $$|\Omega|=\sum_{j=1}^{\infty}\leq 3^n(B+[K]_{H_{q'}})\|f\|_{L^q(\mathbb{R}^n)}^q.$$

For $\text{II}$, use Chebyshev's inequality, the fact that $\int_{Q_j}b_j(y)-(17\sqrt{n})^{\frac{n}{q}}\gamma\alpha\mathbbm{1}_{E_j}(y)\,dy=0$, Fubini's theorem, and H\"older's inequality to estimate
\begin{align*}
\text{II}&\leq 4\alpha^{-1}\int_{\mathbb{R}^n\setminus\Omega}\left|T\left(b-(17\sqrt{n})^{\frac{n}{q}}\gamma\alpha \mathbbm{1}_E\right)(x)\right|dx\\
&\leq 4\alpha^{-1}\sum_{j=1}^{\infty}\int_{\mathbb{R}^n\setminus \Omega}\left|T\left(b_{j}-(17\sqrt{n})^{\frac{n}{q}}\gamma\alpha\mathbbm{1}_{E_{j}}\right)(x)\right|dx\\
&\leq 4\alpha^{-1}\sum_{j=1}^{\infty}\int_{\mathbb{R}^n\setminus  \Omega}\int_{Q_j}|K(x-y)-K(x-c_j)|\left|b_j(y)-(17\sqrt{n})^{\frac{n}{q}}\gamma\alpha\mathbbm{1}_{E_j}(y)\right|dydx\\
&= 4\alpha^{-1}\sum_{j=1}^{\infty}\int_{Q_j}\left(\int_{\mathbb{R}^n\setminus \Omega}|K(x-y)-K(x-c_j)|dx\right)\left|b_j(y)-(17\sqrt{n})^{\frac{n}{q}}\gamma\alpha\mathbbm{1}_{E_j}(y)\right|dy\\
&\leq 4\alpha^{-1}\sum_{j=1}^{\infty}\left\|\int_{\mathbb{R}^n\setminus \Omega}|K(x-y)-K(x-c_j)|dx\right\|_{L^{q'}(Q_j)}\left\|b_j-(17\sqrt{n})^{\frac{n}{q}}\gamma\alpha\mathbbm{1}_{E_j}\right\|_{L^q(\mathbb{R}^n)}\\
&\leq 4\alpha^{-1}\sup_{j\in\mathbb{N}}\left\|\int_{\mathbb{R}^n\setminus\Omega}|K(x-y)-K(x-c_j)|dx\right\|_{L^{q'}\left(Q_j, \frac{dy}{|Q_j|}\right)}\\
&\quad\quad\times\sum_{j=1}^{\infty}|Q_j|^{\frac{1}{q'}}\left\|b_j-(17\sqrt{n})^{\frac{n}{q}}\gamma\alpha\mathbbm{1}_{E_j}\right\|_{L^q(\mathbb{R}^n)}.
\end{align*}
Using the triangle inequality, property (3), and the fact that $|E_j|\leq |Q_j|$, we have $$\left\|b_j-(17\sqrt{n})^{\frac{n}{q}}\gamma\alpha\mathbbm{1}_{E_j}\right\|_{L^q(\mathbb{R}^n)} \leq \|b_j\|_{L^q(\mathbb{R}^n)}+(17\sqrt{n})^{\frac{n}{q}}\gamma\alpha|E_j|^{\frac{1}{q}}\leq 2(17\sqrt{n})^{\frac{n}{q}} \gamma\alpha|Q_j|^{\frac{1}{q}}.$$
Using the above estimate and property (2), we control  
\begin{align*}\hspace{-.2in}
\text{II}\,\le \, &  8(17\sqrt{n})^{\frac{n}{q}}\gamma\sup_{j\in\mathbb{N}}\left\|\int_{\mathbb{R}^n\setminus \Omega}|K(x-y)-K(x-c_j)|dx\right\|_{L^{q'}\left(Q_j, \frac{dy}{|Q_j|}\right)}\sum_{j=1}^{\infty}|Q_j|
\\  \,\le\, & 8(17\sqrt{n})^{\frac{n}{q}}3^n\gamma^{1-q}\alpha^{-q}\|f\|_{L^q(\mathbb{R}^n)}^q\sup_{j\in\mathbb{N}}\left\|\int_{\mathbb{R}^n\setminus \Omega}|K(x-y)-K(x-c_j)|dx\right\|_{L^{q'}\big(Q_j, \frac{dy}{|Q_j|}\big)}.
 \end{align*}

For each $j$, setting $R_j=\frac{\sqrt{n}}{2}l(Q_j)$, we have
$$
Q_j \subseteq  B(c_{j},  R_j) \subseteq  B(c_{j}, 2R_j)  \subseteq  \Omega.
$$ Then the supremum is bounded by 
 $$\sup_{j\in\mathbb{N}}\bigg[ \int_{B(c_{j}, R_j)} 
  \bigg(\int_{\mathbb{R}^n\setminus B(c_{j}, 2R_j)}  |K(x-y)-K(x-c_{j})|   dx
 \bigg)^{q'}\!\! \frac{dy}{|Q_j|} \bigg]^{\frac{1}{q'}},$$
which is bounded by $\left(\frac{\sqrt{n}}{2}\right)^nv_n [K]_{H_{q'}}$ by changing variables $x'=x-c_{j}$, 
$y'=y-c_{j}$ and 
by replacing the supremum over $R_j $ by the supremum over all $R>0$. Therefore $$\text{II}\leq8(17\sqrt{n})^{\frac{n}{q}}\left(\frac{3\sqrt{n}}{2}\right)^nv_n(B+[K]_{H_{q'}})^q\alpha^{-q}\|f\|_{L^q(\mathbb{R}^n)}^q.$$

To control $\text{III}$, use Chebyshev's inequality, the bound of $T$ on $L^s(\mathbb{R}^n)$, the fact that $|E|\leq |\Omega|$, and property (2) to estimate
\begin{align*}
\text{III} &\leq 4^s(17\sqrt{n})^{\frac{ns}{q}}\gamma^s\int_{\mathbb{R}^n}\left|T(\mathbbm{1}_{E})(x)\right|^{s}dx\\
&\leq 4^s(17\sqrt{n})^{\frac{ns}{q}}\gamma^s B^s |E|\\
&\leq 4^s(17\sqrt{n})^{\frac{ns}{q}} |\Omega|\\
&\leq 4^s(17\sqrt{n})^{\frac{ns}{q}}3^n(B+[K]_{H_{q'}})^q\alpha^{-q}\|f\|_{L^q(\mathbb{R}^n)}^q.
\end{align*}

Putting the estimates together, we get $$|\{|Tf|>\alpha\}|\leq \left(2^s+3^n+8(17\sqrt{n})^{\frac{n}{q}}\left(\frac{3\sqrt{n}}{2}\right)^nv_n+4^s(17\sqrt{n})^{\frac{ns}{q}}3^n\right)
\frac{ (B+[K]_{H_{q'}})^q}{\alpha^{q}}\|f\|_{L^q(\mathbb{R}^n)}^q.$$
Since we assumed that $f$ was nonnegative, we must double the constant above to prove the statement for general $f\in L^q(\mathbb{R}^n)$.

When $s=\infty$, set $\gamma=(4(B+[K]_{H_{q'}}))^{-1}$. Then $$\|Tg\|_{L^{\infty}(\mathbb{R}^n)}\leq B\|g\|_{L^{\infty}(\mathbb{R}^n)}\leq B\gamma\alpha\leq \frac{\alpha}{4},$$ so $\left|\left\{|Tg|>\frac{\alpha}{2}\right\}\right|=0$. The part of the argument involving the set $\left\{|Tb|>\frac{\alpha}{2}\right\}$ is the same as in the case $s<\infty$.
\end{proof}

\section{Conclusion}
We end with some remarks and an open question.

\begin{remark}{3}
The conclusions of Theorem 1 and Corollary 1 also follow under the weaker hypothesis that $T$ is bounded from $L^{s,1}(\mathbb{R}^n)$ to $L^{s,\infty}(\mathbb{R}^n)$. Here $L^{s,r}(\mathbb{R}^n)$ is the usual Lorentz space. 
\end{remark}

\begin{remark}{4}
As in the case $q=1$, there are natural vector-valued extensions of Theorem 1 and Corollary 1, in the spirit of \cite{BCP1962}.
\end{remark}

\begin{remark}{5}
Theorem 1 and Corollary 1 are also valid if the original kernel is not of convolution type. In this setting, we say a kernel $K$ defined on $\mathbb{R}^n \times \mathbb{R}^n \setminus \{(x,y):x=y\}$ is in $H_r$ if $$\sup_{R>0} \bigg[\frac{1}{v_nR^n}\int_{|y-y'|\le R} \bigg(
 \int_{|x-y|\ge 2 R}   |K(x,y)-K(x,y')|   \, dx\bigg)^{r} dy \bigg]^{\frac{1}{r}} <\infty,$$ and $$\sup_{R>0} \bigg[\frac{1}{v_nR^n}\int_{|x-x'|\le R} \bigg(
 \int_{|x-y|\ge 2 R}   |K(x,y)-K(x',y)|   \, dx\bigg)^{r} dy \bigg]^{\frac{1}{r}} <\infty,$$ where $v_n$ is the volume of the unit ball $B(0,1)$ in $\mathbb{R}^n$.
\end{remark}

As stated in Remark 2 in the introduction, if $q>1$ and $s < \infty$, then $T$ satisfies strong $L^p(\mathbb{R}^n)$ estimates for $p \in\big(\min (s',q), \max (q',s)  \big)$, and in this case, the interval $\big(\min (s',q), \max (q',s)  \big)$ is properly contained in $(1,\infty)$. 

Let $q>1$ and $s< \infty$. 
As of this writing, we are unable to establish whether  the   interval
  $\big(\min (s',q), \max (q',s)  \big)$ is the largest interval $(a,b)$ for which 
 an operator $T$ with kernel in $H_{q'}$ is bounded on $L^p(\mathbb{R}^n)$ for all $p \in (a,b)$. This certainly relates to the existence of examples of kernels in 
 $H_{q_1}$ but not in $H_{q_2}$ for $q_1<q_2$.


\end{document}